\documentclass[11 pt,reqno,twoside]{article}
\usepackage{graphicx}
          \DeclareGraphicsExtensions{.pdf}
\usepackage[utf8]{inputenc}
\pdfoutput=1

\usepackage{enumerate}
\usepackage{amsmath}
\usepackage{latexsym}
\usepackage{amssymb}

\newcommand{\rank}{\mathop{\rm rank}\nolimits}

\def\Box{\rule{2mm}{2mm}}

\font\tenbbb=msbm10 scaled\magstephalf \font\sevenbbb=msbm7
\font\fivebbb=msbm5
\newfam\bbbfam
\textfont\bbbfam=\tenbbb \scriptfont\bbbfam=\sevenbbb
\scriptscriptfont\bbbfam=\fivebbb

\newtheorem{thm}{Theorem}
\newtheorem{cor}[thm]{Corollary}
\newtheorem{lem}[thm]{Lemma}
\newtheorem{prop}[thm]{Proposition}
\newtheorem{defn}{Definition}

\newtheorem{ex}{Example}

\def\p{\begin{exercici}}
\def\fp{\end{exercici}}

\newcommand{\kera}{\mathop{\rm Ker}\nolimits}

\newcommand{\tr}{\mathop{\rm tr}\nolimits}
\newcommand{\sign}{\mathop{\rm sign}\nolimits}

\newcommand{\be}{\begin{exercici}}
\newcommand{\ee}{\end{exercici}}

\newcommand{\fidibuix}{\end{minipage}}

\font\5=cmbx10 scaled \magstep 3 \font\4=cmbx10 scaled \magstep 2
\title{
\sc \bf Unobservable Planar Bimodal Linear Systems: Miniversal
Deformations, Controllability and
Stabilization\thanks{Partially supported by DGICYT
  MTM2010-19356-C02-02.}
}
\author{\sc J.~Ferrer, M.~D.~Magret, J.~R.~Pacha, M.~Pe\~ na\\\\
Departament de Matem\`{a}tica Aplicada I.\\
E.T.S. Enginyeria Industrial de Barcelona. UPC\\
Diagonal 647. 08028 Barcelona. Spain\\
josep.ferrer@upc.edu, m.dolors.magret@upc.edu,\\
juan.ramon.pacha@upc.edu, marta.penya@upc.edu}
\date{}

\font\pt=cmr9 
\begin{document}
\maketitle
\begin{abstract}
We consider the set of bimodal linear systems consisting of two linear
dynamics acting on each side of a given hyperplane, assuming continuity
along the separating hyperplane. Focusing on the unobservable planar
ones, we obtain a simple explicit characterization of
controllability. Moreover, we apply the canonical forms of these systems
depending on two state variables to obtain explicitly miniversal
deformations, to illustrate bifurcation diagrams and to prove that the
unobservable controllable systems are stabilizable.
\end{abstract}

%

\section{Introduction}\label{sec:intro}
Piecewise linear systems have attracted the interest of the researchers
in recent years by their wide range of applications, as well as by the
possible theoretical approaches. In both directions, the non-generic
case of unobservable systems presents interesting particularities. They
appear in a natural way in parameterized families, bifurcations,..., and
special theoretical tools are needed to study their properties.

In this paper, we consider the set of bimodal linear systems consisting
of two linear dynamics acting on each side of a given hyperplane,
assuming continuity along the separating hyperplane. In the space of
triples of matrices defining those systems, we consider the natural
equivalence relation defined by changes of bases in the state space
which preserves the hyperplanes parallel to the separating
hyperplane. The fact that this equivalence relation can be viewed as the
one induced by the action of a Lie group allows to apply Arnold's
techniques concerning versal deformations, stratification, etc.

Canonical forms, representative for each equivalence class, are a basic
tool in order to simplify the computations, because instead of the given
matrices their canonical forms can be used. More concretely, we present
an explicit expression of the canonical form of a given triple of matrices
when the number of state variables is two or three, which are the most
commonly found in applications (see for example~\cite{CarmonaFPT02},
\cite{CarmonaFPT06}, \cite{DiBernardoPP08}, \cite{DiBernardoBCK08},
\cite{FreirePR05}). 
We use these canonical forms to obtain explicitly miniversal
deformations and to illustrate local bifurcation diagrams always
focusing on the unobservable case.

Secondly, we obtain a simple explicit characterization of the
controllability of a planar bimodal system, starting on the implicit one
in~\cite{CamlibelHS}. Indeed, we prove that for $n=2$ one of the
conditions there can be avoided (Corollary~\ref{cor:contro3}), and the
other one can be reformulated in a very simple way
(Proposition~\ref{prop:contro2}).

Moreover, in~\cite{CamlibelHS} one asks if controllable bimodal linear
systems can be stabilized by means of a feedback (the same for both
subsystems), generalizing the well-known result for a single
system. Here the reformulation of controllability
(Corollary~\ref{cor:contro4}) and the above canonical forms are used to
prove that this is true in the unobservable planar case
(Theorem~\ref{teor:stabi}).
%

The structure of the paper is as follows. In Section~\ref{sec:BPWLS} we
provide the canonical forms obtained for a bimodal linear system in the
cases where the number of state variables is two, and recall the
dimensions of the orbits and the strata. As an application, in
Section~\ref{sec:miniversal} we compute miniversal deformations and
apply them to obtain local bifurcation diagrams. Section~\ref{sec:contr}
is devoted to obtain a simple explicit characterization of the
controllability and to prove that then the system is stabilizable by
feedback.

Throughout the paper, $\mathbb{R}$ will denote the set of real numbers,
$M_{n\times m}(\mathbb{R})$ the set of matrices having $n$ rows and $m$
columns and entries in $\mathbb{R}$ (in the case where $n=m$, we will
simply write $M_n(\mathbb{R})$) and  $Gl_n(\mathbb{R})$ the group of
non-singular matrices in $M_n(\mathbb{R})$. Finally, we will denote by
$e_1,\dots ,e_n$ the natural basis of the Euclidean space
${\mathbb{R}}^n$.

\section{Canonical Forms}\label{sec:BPWLS}
Let us consider a bimodal linear dynamical system given by
\begin{displaymath}
\begin{cases}
   \dot{\boldsymbol{x}}(t) = A_{1}\boldsymbol{x}(t) + B_{1},\\
         \boldsymbol{y}(t) = C\boldsymbol{x}(t),
\end{cases} \
        \text{ if }\;
        \boldsymbol{y}(t)\le 0,\quad
\begin{cases}
  \dot{\boldsymbol{x}}(t) = A_{2} \boldsymbol{x}(t) + B_{2},\\
        \boldsymbol{y}(t) = C\boldsymbol{x}(t),
\end{cases}\ \text{ if }\;
        \boldsymbol{y}(t)\ge 0
\end{displaymath}
where $A_{1},A_{2}\in M_{n}(\mathbb{R})$; $B_{1},B_{2}\in M_{n\times 1}
(\mathbb{R})$; $C\in M_{1\times n}(\mathbb{R})$. Let us
assume that the dynamics is continuous along the separating hyperplane
$H=\{\boldsymbol{x}\in\mathbb{R}^{n}: C\boldsymbol{x}=0\}$. For
simplicity, we will consider that $C=(1\;0\ldots 0)\in
M_{1\times n}({{\mathbb{R}}})$. Hence  $H =\{\boldsymbol{x}\in
\mathbb{R}^{n}:x_{1}=0\}$ and continuity along $H$ is equivalent
to:
\begin{displaymath}
   B_{2}= B_{1},\qquad A_{2}e_{i}= A_{1}e_{i},\quad 2\le i\le n.
\end{displaymath}
We will write from now on $B=B_1=B_2$.

\begin{defn}
In the above conditions, we say that the triple of matrices
$(A_{1},A_{2},B)$ defines a bimodal piecewise linear system. Throughout
the paper, $\mathcal{X}$ will denote the set of these triples
\begin{displaymath}
  \mathcal{X}=\left\{(A_{1},A_{2},B)\in M_{n}(\mathbb{R})\times
    M_{n}(\mathbb{R})\times M_{n\times 1}(\mathbb{R})
    \left\vert\; A_{2} e_{i} =A_{1}e_{i},\: 2\leq i\leq n\right.\right\}
\end{displaymath}
which is obviously a ($n^2+2n$)-differentiable manifold.

The system is called \emph{observable} if
\begin{displaymath}
   \rank \left (\begin{array}{c}C\\CA_{i}\\\cdots\\CA_{i}^{n-1}
\end{array}\right)=n,\quad i=1,2.
\end{displaymath}
\end{defn}

The basis changes in the state variables space preserving the
hyperplanes $x_{1}(t)=k$ will be called \emph{admissible basis
  changes}. Thus, they are basis changes given by a matrix $S\in
Gl_{n}(\mathbb{R})$,
\begin{displaymath}
   S=\left (\begin{array}{cc}
               1 & 0\\
               U & T
            \end{array}\right),\quad T\in Gl_{n-1}(\mathbb{R}).
 \end{displaymath}
We consider the equivalence relation in the set of triples of matrices
$\mathcal{X}$ which corresponds to admissible basis changes:

\begin{defn}
We write
\begin{displaymath}
   \mathcal{S}:=\left\{S\in\, Gl_{n}(\mathbb{R})\left\vert\,
   S=\left(\begin{array}{cc}
              1 & 0\\
              U & T
           \end{array}\right),\; T\in Gl_{n-1}({\mathbb{R}})\right.
           \right\}.
\end{displaymath}
\end{defn}
Then, $(A_{1}, A_{2}, B),(A'_{1}, A'_{2}, B')\in \mathcal{X}$ are said
to be \emph{equivalent} if there exists a matrix $S\in {\cal S}$
(representing an admissible basis change) such that $
(A'_{1},A'_{2},B')=(S^{-1}A_{1}S,S^{-1}A_{2}S,S^{-1}B)$.

Notice that the matrix $C$ is not involved in this definition since
$CS=C$ for any $S\in {\cal S}$.

We remark that equivalence classes are actually the orbits with regard
to the action of the Lie group ${\cal S}$ on the differentiable manifold
${\cal X}$,
\begin{displaymath}
\begin{array}{rl}
   \alpha: {\cal S}\times {\cal X} & \longrightarrow {\cal X}
\end{array}
\end{displaymath}
defined by
\begin{displaymath}
   \alpha(S, {\cal X})=(S^{-1}A_{1}S,S^{-1}A_{2}S,S^{-1}B)
\end{displaymath}
Given any triple of matrices $(A_{1}, A_{2}, B)\in {\cal X}$, we will
denote by ${\cal O}(A_{1},A_{2},B)$ its orbit (or equivalence class).

As an application of the Closed Orbit Lemma (see~\cite{Hum}), we deduce
that equivalence classes are differentiable manifolds. Namely,
any equivalence class is a locally closed differentiable submanifold of
${\cal X}$ and its boundary is a union of equivalence classes or orbits
of strictly lower dimension. In particular, equivalence classes or
orbits of minimal dimension are closed.

In this section we summarize the results in the previous works~\cite{MP}
and~\cite{MPP} which will be used in the sequel. The first goal is
obtaining for each triple $(A_{1},A_{2},B)$ a canonical reduced form which
characterizes its equivalence class. In
\texttt{http://www.ma1.upc.edu/$\sim$joanr/html/cfbpwls.html} the reader
can find a MAPLE program which allows to obtain the canonical form of a
triple $(A_{1}, A_{2}, B)$ for the cases $n=2$ and $n=3$, which are the
most commonly systems found in practice. Furthermore, one obtains an
admissible basis change $S\in {\cal S}$ which transforms the initial
triple $(A_{1}, A_{2}, B)$ into its canonical form.

When listing canonical forms, it is necessary that the coefficients
appearing in them as well as the conditions used to distinguish the
different types do not depend on the admissible basis which one
considers, that is to say, they are preserved under admissible basis
changes $S\in {\cal S}$. It is well-known that $\tr A_{1}$, $\tr A_{2}$,
$\det A_{1}$, $\det A_{2}$ are invariant under any basis change
$S\in Gl_{n}(\mathbb{R})$. We focus on the additional invariants when
only admissible basis changes $S\in \mathcal{S}$ are considered.

\begin{defn}
A real number (respectively a property) associated to a triple
$(A_{1}, A_{2}, B)$ is called \emph{$\mathcal{S}$-invariant} if it is
preserved by admissible basis changes, that is to say, it has the same
value (respectively it is also true) for any other triple
$(A'_{1}, A'_{2},B')$ $\mathcal{S}$-equivalent to the given one.
\end{defn}
For example, it is obvious that:
\begin{prop}\label{prop:pres}
 They are $\mathcal{S}$-invariant:
\begin{enumerate}[(i)]
\item The top coefficient $b_{1}$ in $B$ ($b_{1}=CB$).
\item\label{item:matrixC}  The matrix $C$.
\item  The condition of $(A_1,A_2,B)$ being observable.
\end{enumerate}
\end{prop}
Next proposition details the remainder $\mathcal{S}$-invariants that we
will use for $n=2$. In order to that, we define:

\begin{defn}
Given a triple
\begin{displaymath}
   A_{1}=\left(\begin{array}{cc}
                  a_{1} & a_{3}\\
                  a_{2} & a_{4}
               \end{array} \right),
   A_{2}=\left(\begin{array}{cc}
                  \gamma_{1} & a_{3}\\
                  \gamma_2&a_4\end{array} \right),
       B=\left(\begin{array}{c}
                   b_{1}\\
                   b_{2}\end{array}\right)
\end{displaymath}
we write:
\begin{displaymath}
\begin{array}{rl}
     & \Delta_{0} = \det\left(\begin{array}{cc}
                                 a_{3} & b_{1}\\
                                 a_{4} & b_{2}
                              \end{array} \right)
                  = a_{3} b_{2} - a_{4} b_{1}\\
     & \Delta_{12} = a_{2} (a_{4}-\gamma_{1})-\gamma_{2} (a_{4}-a_{1})\\
     & \Delta_{1} = b_{1} a_{2} + (a_{4} - a_{1}) b_{2}\\
     & \Delta_{2} = b_{1}\gamma_{2} + (a_{4} - \gamma_{1}) b_{2}
\end{array}
\end{displaymath}
\end{defn}

\begin{lem}\label{lem:lemn1}
The above triple is \emph{unobservable} if and only if $a_{3}=0$. In
this case we have:

\begin{enumerate}[(1)]
\item\label{item:det} $\det \left(\left(\begin{array}{c}
                                a_{1}-\gamma_{1}\\
                                a_{2}-\gamma_{2}
                        \end{array}\right)\left| A_{i}
                        \left(\begin{array}{c}
                                 a_{1} - \gamma_{1}\\
                                 a_{2} - \gamma_{2}
                               \end{array}\right)\right.\right)
                        = (a_{1}-\gamma_{1})\Delta_{12}$,\quad $i=1,2$\\
      $\det\left(\begin{array}{c|c}
                    B & A_{i}B
                 \end{array}\right) = b_{1}\Delta_{i}$,\quad $i=1,2$
\item The action of $S\in\mathcal{S}$ transforms $\Delta_{1}$,
  $\Delta_{2}$, $\Delta_{12}$ respectively into:
\begin{displaymath}
\frac{1}{\det S}\Delta_{1},\quad\frac{1}{\det S}\Delta_{2},\quad\frac{1}{\det
  S}\Delta_{12}
\end{displaymath}
In particular, it is $\mathcal{S}$-invariant the sign (positive,
negative or zero):
\begin{displaymath}
\sign(\Delta_{1}\Delta_{2})
\end{displaymath}
\end{enumerate}
\end{lem}

\noindent {\bf Proof.}
Clearly,
\begin{equation}
   \left(\begin{array}{c}
            C\\
            CA_{1}
         \end{array}\right) =
   \left(\begin{array}{cc}
            1 & 0\\
            a_{1} & a_{3}
          \end{array}\right),\quad
   \left(\begin{array}{c}
            C\\
            CA_{2}
         \end{array}\right) =
   \left(\begin{array}{cc}
            1 & 0\\
            \gamma_{1} & a_{3}
         \end{array}\right)\label{eq:matCCA}
\end{equation}
do not have maximal rank when $a_{3}=0$. Then:
\begin{enumerate}[(1)]
\item It is a straightforward computation.
\item If $a_3=0$, then $a_1$ and $\gamma_1$ are eigenvalues of $A_1$ and $A_2$.\\
The action of $S$ transforms the matrices in~\eqref{eq:matCCA} into
their left product by $S^{-1}$.
\end{enumerate}
\hfill $\Box$

\begin{prop}
With the above notation, the following table summarizes some
$\mathcal{S}$-invariant numbers and properties, as well as the
hypotheses for each one:

\bigskip
\centerline{ \begin{tabular}{lccc}
 &   & \underline{numbers} & \underline{properties} \\
(1)  &                     & $\Delta_{0}$ & $a_{3} = 0$\\
(2)  & $b_{1} = 0$          &             & $b_{2} = 0$\\
(3)  & $a_{3} = 0$          & $a_{1},\gamma_{1},a_{4}$ &
                             $\Delta_{12} = 0, \Delta_{1} = 0,
                              \Delta_{2} =0$\\
(3') & $a_{3} = 0, \Delta_{12} \neq 0$ & $\Delta_{1}/\Delta_{12},
                                        \Delta_{2}/\Delta_{12}$ & \\
(4)  & $a_{3} = 0, a_{1}=a_{4}$         &                        &
                                                      $a_{2} = 0$\\
(4') & $a_{3} = 0,\gamma_{1}=a_{4}$     &         & $\gamma_{2} = 0$\\
(5)  &  $b_{1}=0, a_{3} = 0, a_{1} = a_{4}$        & $b_{2}/a_{2}$  & \\
(5') & $b_{1}=0, a_{3} = 0, \gamma_{1} =a_{4}$ & $b_{2}/\gamma_{2}$ & \\
\end{tabular} }
\end{prop}

\noindent
{\bf Proof.}
\begin{enumerate}[(1)]
\item\label{item:Delta0a30} The $\mathcal{S}$-action on $A_{1}$ and
  $B$ can be formulated as:
\begin{displaymath}
S^{-1}(A_{1},B)\left(\begin{array}{c|c}
                      S & 0\\ \hline
                      0 & 1\end{array}\right).
\end{displaymath}
For $n=2$, one has
\begin{displaymath}
\left(\begin{array}{cc}
         1 & 0\\
         u & t
      \end{array}\right)^{-1}
\left(\begin{array}{ccc}
         a_{1} & a_{3} & b_{1}\\
         a_{2} & a_{4} & b_{2}\end{array}\right)
\left(\begin{array}{ccc}
         1 & 0 & 0\\
         u & t & 0\\
         0 & 0 & 1\end{array}\right) =
\left(\left.\begin{array}{c}{\ast}\\ {\ast}\end{array}\right|
            \left(\begin{array}{cc}
                     1 & 0\\
                     u & t\end{array}\right)^{-1}
            \left(\begin{array}{cc}
                     a_{3} & b_{1}\\
                     a_{4} & b_{2}\end{array}\right)
            \left(\begin{array}{cc}
                     t & 0\\
                     0 & 1\end{array}\right)\right ).
\end{displaymath}
Therefore, $\Delta_{0}$ is $\mathcal{S}$-invariant:
\begin{displaymath}
\det\left(\left(\begin{array}{cc}
                   1 & 0\\
                   u & t\end{array}\right)^{-1}
          \left(\begin{array}{cc}
                   a_{3} & b_{1}\\
                   a_{4} & b_{2}\end{array}\right)
          \left(\begin{array}{cc}
                   t & 0\\
                   0 & 1\end{array}\right)\right) =
\det\left(\begin{array}{cc}
                   a_{3} & b_{1}\\
                   a_{4} & b_{2}\end{array}\right).
\end{displaymath}
We have seen that $a_{3}\neq 0$ if and only if
\begin{displaymath}
\mbox{\rm rk}\left(\begin{array}{c}
                      C\\
                      CA_{i}\end{array}\right) = 2, \quad i=1,2,
\end{displaymath}
which is $\mathcal{S}$-invariant.
\item If $b_{1}=0$, then
\begin{displaymath}
S^{-1}\left(\begin{array}{c}
               0\\
               b_{2}\end{array}\right) =
\frac{1}{t}\left(\begin{array}{rc}
                    t & 0\\
                   -u & 1\end{array}\right)
            \left(\begin{array}{c}
                    0\\
                    b_{2}\end{array}\right) =
\frac{1}{t}\left(\begin{array}{c}
                    0\\
                    b_{2}\end{array}\right).
\end{displaymath}
\item\label{item:a3eq0} If $a_{3}=0$, then $a_{1}, a_{4}, \gamma_{1}$
  are the eigenvalues of $A_{1}$ and $A_{2}$. Then:
\begin{displaymath}
   \kera(A_{1} - a_{1}I)=\kera(A_{2} - \gamma_{1} I)
\end{displaymath}
  if and only if
\begin{displaymath}
  \mbox{\rm rk}\left(\begin{array}{cc}
                        a_{2} & a_{4} - a_{1}\\
                        \gamma_{2} & a_{4} - \gamma_{1}
                     \end{array}\right) = 1
\end{displaymath}
or, equivalently,
\begin{displaymath}
   0 = \det\left(\begin{array}{cc}
                    a_{2} & a_{4} - a_{1}\\
                   \gamma_{2} & a_{4} - \gamma_{1}
                 \end{array}\right) = \Delta_{12}.
\end{displaymath}
In a similar way, $\Delta_{1} = 0$ if and only if
\begin{displaymath}
   \left(\begin{array}{c}
            b_{1}\\
            b_{2}\end{array}\right)\in\kera(A_{1} - a_{1} I)
\end{displaymath}
and $\Delta_{2} = 0$ if and only if
\begin{displaymath}
   \left(\begin{array}{c}
            b_{1}\\
            b_{2}\end{array}\right)\in\kera(A_{2}-\gamma_{1}I).
\end{displaymath}
\item[(3')] It follows from~(\ref{item:a3eq0}) and the above lemma.
\item[(4),](4') Clearly, if $a_{3}=0$ and $a_{1}=a_{4}$, then $a_{2}=0$
  if and only if $A_{1}$ diagonalizes. And analogously for
  $\gamma_{2}=0$.
\item[(5),](5') Returning to the formulation in~(\ref{item:Delta0a30}):
\begin{displaymath}
\left(\begin{array}{cc}
         1 & 0\\
         u & t\end{array}\right)^{-1}
\left(\begin{array}{ccc}
         a_{1} & 0     & 0\\
         a_{2} & a_{1} & b_{2}\end{array}\right)
\left(\begin{array}{ccc}
         1 & 0 & 0\\
         u & t & 0\\
         0 & 0 &1\end{array}\right) =
\left(\begin{array}{ccc}
         a_{1}   & 0     & 0\\
         a_{2}/t & a_{1}  & b_{2}/t\end{array}\right)
\end{displaymath}
and analogously for $\gamma_{1}=a_{4}$.
 \end{enumerate}
\hfill $\Box$

By means of the above $\mathcal{S}$-invariants, one may list the
possible canonical forms and the classification criteria:
\begin{enumerate}\itemindent=10pt
\item[(CF1)$\hphantom{a}$] $a_{3}\neq 0$
\item[{}] $\left(\begin{array}{cc}
                    \tr A_{1}  & 1\\
                    \det A_{1} & 0\end{array}\right),
           \left(\begin{array}{cc}
                    \tr A_{2}  & 1\\
                    \det A_{2} & 0 \end{array}\right),
           \left(\begin{array}{c}
                    b_{1}\\
                    \Delta_{0}\end{array}\right)$
\item[(CF2)$\hphantom{a}$] $a_{3}=0, a_{1}\neq a_{4}, \gamma_{1}\neq
                            a_{4}, \Delta_{12}\neq 0$
\item[{}]  $\left(\begin{array}{cc}
                     a_{1} & 0\\
                        0 & a_{4}\end{array}\right),
            \left(\begin{array}{cc}
                    \gamma_{1} & 0\\
                         1    & a_{4}\end{array}\right),
            \left(\begin{array}{c}
                    b_{1}\\
                    -\Delta_1/\Delta_{12}\end{array}\right)$
\newpage
\item[(CF3)$\hphantom{a}$] $a_{3}=0, a_{1}\neq a_{4},
                           \gamma_{1}\neq a_{4}, \Delta_{12}=0,
                           \Delta_{1}\neq 0$
\item[{}] $\left(\begin{array}{cc}
                    a_{1} & 0\\
                       0 & a_{4}\end{array}\right),
           \left(\begin{array}{cc}
                   \gamma_{1} & 0\\
                    0        & a_{4}\end{array}\right),
           \left(\begin{array}{c}
                    b_{1}\\
                    1\end{array}\right)$
\item[(CF4)$\hphantom{a}$] $a_{3}=0, a_{1}\neq a_{4},
                           \gamma_{1}\neq a_{4}, \Delta_{12}=0,
                           \Delta_{1}=0$
\item[{}] $\left(\begin{array}{cc}
                    a_{1} & 0\\
                       0 & a_{4}\end{array}\right),
           \left(\begin{array}{cc}
                   \gamma_{1} & 0\\
                       0     & a_{4}\end{array}\right),
           \left(\begin{array}{c}
                       b_{1}\\
                       0\end{array}\right)$
\item[(CF5)$\hphantom{a}$] $a_{3} = 0, a_{1} = a_{4},
                           \gamma_{1}\neq a_{4}, a_{2}\neq 0$
\item[{}]  $\left(\begin{array}{cc}
                     a_{4} & 0\\
                        1 & a_{4}\end{array}\right),
            \left(\begin{array}{cc}
                     \gamma_{1} & 0\\
                             0 & a_{4}\end{array}\right),
            \left(\begin{array}{c}
                     b_{1}\\
                     \Delta_{2}/\Delta_{12}\end{array}\right)$
\item[(CF5')] $a_{3} = 0, a_{1}\neq a_{4}, \gamma_{1} = a_{4},
               \gamma_{2}\neq 0$
\item[{}] $\left(\begin{array}{cc}
                    a_{1} & 0\\
                       0 & a_{4}\end{array}\right),
           \left(\begin{array}{cc}
                    a_{4} & 0\\
                       1 & a_{4}\end{array}\right ),
           \left(\begin{array}{c}
                    b_{1}\\
                  -\Delta_{1}/\Delta_{12}\end{array}\right)$
\item[(CF6)$\hphantom{a}$] $a_{3} = 0, a_{1} = a_4,
                           \gamma_{1}\neq a_{4}, a_{2} = 0,
                           \Delta_{2}\neq 0$
\item[{}] $\left(\begin{array}{cc}
                    a_{4} & 0\\
                       0 & a_{4}\end{array}\right),
           \left(\begin{array}{cc}
                    \gamma_{1} & 0\\
                            0 & a_{4}\end{array}\right),
           \left(\begin{array}{c}
                    b_{1}\\
                    1\end{array}\right)$
\item[(CF6')] $a_{3} = 0, a_{1}\neq a_{4}, \gamma_{1} = a_{4},
               \gamma_{2} = 0, \Delta_{1}\neq 0$
\item[{}] $\left(\begin{array}{cc}
                    a_{1} & 0\\
                       0 & a_{4}\end{array}\right),
           \left(\begin{array}{cc}
                    a_{4} & 0\\
                       0 & a_{4}\end{array}\right),
           \left(\begin{array}{c}
                    b_{1}\\1
           \end{array}\right)$
\item[(CF7)$\hphantom{a}$] $a_{3} = 0, a_{1} = a_{4},
                           \gamma_{1}\neq a_{4},
                           a_{2} = 0, \Delta_{2} = 0$
\item[{}]  $\left(\begin{array}{cc}
                     a_{4} & 0\\
                        0 & a_{4}\end{array}\right),
            \left(\begin{array}{cc}
                     \gamma_{1} & 0\\
                             0 & a_{4}\end{array}\right),
            \left(\begin{array}{c}
                     b_{1}\\0\end{array}\right)$
\item[(CF7')] $a_{3} = 0, a_{1}\neq a_{4}, \gamma_{1} = a_{4},
              \gamma_{2} = 0, \Delta_{1} = 0$
\item[{}] $\left(\begin{array}{cc}
                    a_{1} & 0\\
                       0 & a_{4}\end{array}\right),
           \left(\begin{array}{cc}
                    a_{4} & 0\\
                       0 & a_{4}\end{array}\right),
           \left(\begin{array}{c}
                    b_{1}\\
                      0\end{array}\right)$
\item[(CF8)$\hphantom{a}$] $a_{3} = 0, a_{1} = a_{4} = \gamma_{1},
                            a_{2}\neq 0, \gamma_{2}\neq 0,
                            b_{1}\neq 0$
\item[{}] $\left(\begin{array}{cc}
                    a_{4} & 0\\
                       1 & a_{4}\end{array}\right),
           \left(\begin{array}{cc}
                    a_{4} & 0\\
                    \gamma_{2}/a_{2} & a_{4}\end{array}\right),
           \left(\begin{array}{c}
                    b_{1}\\0 \end{array}\right)$
\item[(CF9)$\hphantom{a}$] $a_{3} = 0, a_{1} = a_{4} = \gamma_{1},
                            a_{2}\neq 0, \gamma_{2}\neq 0, b_{1} = 0$
\item[{}] $\left(\begin{array}{cc}
                    a_{4} & 0\\
                       1 & a_{4}\end{array}\right),
           \left(\begin{array}{cc}
                    a_{4} & 0\\
                    \gamma_{2}/a_{2} & a_{4}\end{array}\right),
           \left(\begin{array}{c}
                    0\\
                    b_{2}/a_{2}\end{array}\right)$
\item[(CF10)$\hphantom{a}$] $a_{3} = 0, a_{1} = a_{4} = \gamma_{1},
                             a_{2}\neq 0, \gamma_{2} = 0, b_{1}\neq 0$
\item[{}] $\left(\begin{array}{cc}
                    a_{4} & 0\\
                       1 & a_{4}\end{array}\right),
           \left(\begin{array}{cc}
                    a_{4} & 0\\
                       0 & a_{4}\end{array}\right),
           \left(\begin{array}{c}
                    b_{1}\\0
           \end{array}\right)$
\newpage
\item[(CF10')] $a_{3} = 0, a_{1} = a_{4} = \gamma_{1}, a_{2}=0,
               \gamma_{2}\neq 0, b_{1}\neq 0$
\item[{}] $\left(\begin{array}{cc}
                    a_{4} & 0\\
                       0 & a_{4}\end{array}\right),
           \left(\begin{array}{cc}
                    a_{4} & 0\\
                       1 & a_{4}\end{array}\right),
           \left(\begin{array}{c}
                    b_{1}\\
                       0\end{array}\right)$
\item[(CF11)$\hphantom{a}$] $a_{3} = 0, a_{1} = a_{4} = \gamma_{1},
                             a_{2}\neq 0, \gamma_{2} = 0, b_{1} = 0$
\item[{}] $\left(\begin{array}{cc}
                   a_{4} & 0\\
                      1 & a_{4}\end{array}\right),
           \left(\begin{array}{cc}
                   a_{4} & 0\\
                     0  & a_{4}\end{array}\right),
           \left(\begin{array}{c}
                   0\\
                   b_{2}/a_{2}\end{array}\right)$
\item[(CF11')] $a_{3} = 0, a_{1} = a_{4} = \gamma_{1}, a_{2}=0,
                \gamma_{2}\neq 0, b_{1} = 0$
\item[{}] $\left(\begin{array}{cc}
                    a_{4} & 0\\
                       0 & a_{4}\end{array}\right ),
           \left(\begin{array}{cc}
                    a_{4} & 0\\
                       1 & a_{4}\end{array}\right ),
           \left(\begin{array}{c}
                    0\\
                    b_{2}/\gamma_{2}\end{array}\right)$
\item[(CF12)$\hphantom{a}$] $a_{3} = 0, a_{1} = a_{4} = \gamma_{1},
                             a_{2} = 0, \gamma_{2} = 0, b_{1}\neq 0$
\item[{}] $\left(\begin{array}{cc}
                    a_{4} & 0\\
                       0 & a_{4}\end{array}\right),
           \left(\begin{array}{cc}
                    a_{4} & 0\\
                       0 & a_{4}\end{array}\right),
           \left(\begin{array}{c}
                    b_{1}\\
                       0\end{array}\right)$
\item[(CF13)$\hphantom{a}$] $a_{3} = 0, a_{1} = a_{4} = \gamma_{1},
                             a_{2}=0, \gamma_{2} = 0, b_{1} = 0,
                             b_{2}\neq 0$
\item[{}] $\left(\begin{array}{cc}
                    a_{4} & 0\\
                       0 & a_{4}\end{array}\right),
           \left(\begin{array}{cc}
                    a_{4} & 0\\
                       0 & a_{4}\end{array}\right),
          \left(\begin{array}{c}
                    0\\
                    1\end{array}\right)$
\item[(CF14)$\hphantom{a}$] $a_{3} = 0, a_{1} = a_{4} = \gamma_{1},
                             a_{2} = 0, \gamma_{2} = 0, b_{1} = 0,
                             b_{2} = 0$
\item[{}] $\left(\begin{array}{cc}
                    a_{4} & 0\\
                       0 & a_{4}\end{array}\right),
           \left(\begin{array}{cc}
                    a_{4} & 0\\
                       0 & a_{4}\end{array}\right),
           \left(\begin{array}{c}
                    0\\
                    0\end{array}\right)$
\end{enumerate}

Finally, we list the dimension of the orbits and the strata for each
case. We recall that each stratum is the union of the orbits of the
same type when the parameters appearing in the canonical form
vary. In~\cite{MPP} one proves that these sets are differentiable
manifolds.

\bigskip

\centerline{ \begin{tabular}{|c|c|c|}
  \hline
  Canonical form & Dimension of the orbit& Dimension of the stratum \\
  \hline
  CF1 & 2 & 8\\\hline
  CF2 & 2 & 7 \\\hline
  CF3 & 2 & 6\\\hline
 CF4 & 1 & 5 \\\hline
CF5, CF5' & 2 & 6 \\\hline
 CF6, CF6' & 2 & 5 \\\hline
  CF7, CF7' & 1 & 4 \\\hline
  CF8 & 2 & 5 \\\hline
  CF9 &  1 & 4 \\\hline
  CF10, CF10' & 2 & 4 \\\hline
CF11, CF11' &  1 & 3 \\\hline
CF12 & 1 & 3\\\hline
CF13 & 1 & 1 \\\hline
CF14 & 0 & 1 \\ \hline
\end{tabular} }

\section{Miniversal deformations and bifurcation
  diagrams}\label{sec:miniversal}
Versal deformations provide all possible structures which arise when
small perturbations act and can be applied to the study of
singularities and bifurcations. Here we will use them in order to
detail the stratification of the unobservable perturbations of a
given triple. The main definitions and results about deformations and versality
can be found in \cite{Arnold71} and \cite{tan}. Here we re-write them down,
adapted to our particular case.

\begin{defn} A \emph{deformation} of $(A_{1}, A_{2}, B)\in {\cal X}$
  is a differentiable map $\varphi:{U} \longrightarrow {\cal X}$, with
  $U$ an open neighbourhood of the origin ${\mathbb{R}}^d$, such that
$\varphi(0)=(A_1,A_2,B)$.

A deformation $\varphi:{U}\longrightarrow {\cal X}$ of $(A_1,A_2,B)$ is
called  \emph{versal} at $0$ if for any other deformation of
$(A_{1}, A_{2}, B)$, $\psi: {V} \longrightarrow {\cal X}$, there
exists a neighbourhood ${V}' \subseteq {V}$ with $0\in
{V}'$, a differentiable map $\gamma: {V'} \longrightarrow {U}$
with $\gamma(0)=0$ and a deformation of the identity $I \in
{\cal S}$, $\theta:{V}' \longrightarrow {\cal S}$,
such that $\psi(\mu)= \alpha(\theta(\mu),\varphi(\gamma(\mu)))$
for all $\mu \in {V}'$.

A versal deformation with minimal number of parameters $d$ is
called  \emph{miniversal} deformation.
\end{defn}

From the description of the normal space below, the dimension of
miniversal deformations may be computed. Even more, a miniversal
deformation can be obtained from a basis of the normal space to the
orbit of a given triple. This miniversal deformation is usually called
{\it orthogonal} miniversal deformation.

\begin{thm}{Let us denote by
$N_{(A_{1}, A_{2}, B)}{\cal O}(A_{1}, A_{2}, B)$ the normal
space to the orbit of the triple $(A_{1}, A_{2}, B)$ at
$(A_{1}, A_{2}, B)$ with regard to some scalar product in
${\cal X}$. Then, the mapping
\begin{displaymath}
\kern-.2truecm
\begin{array}{rl}
{\mathbb{R}}^d & \kern-.3truecm \longrightarrow {\cal X}\\
(\eta_{1},\dots, \eta_d) & \kern-.3truecm \longrightarrow
(A_{1}, A_{2}, B)+ \eta_{1} V_{1} +\dots +\eta_{d} V_{d}
\end{array}
\end{displaymath}
where $\{V_{1},\dots ,V_{d}\}$ is any basis of the vector space
$N_{(A_{1}, A_{2}, B)}{\cal O}(A_{1}, A_{2}, B)$  is a miniversal deformation of
$(A_{1}, A_{2}, B)$.}
\end{thm}

In \cite{MPP} the authors provided a description of the linear equations
system which leads to a way for computing a basis of the normal space.

\begin{prop} (\cite{MPP})
We consider the following scalar product in ${\cal X}:$
\begin{displaymath}
  \langle (A_{1}, A_{2}, B), (A'_{1}, A'_{2}, B')\rangle =
  \mbox{\rm tr}\,(A_{1}^{t} A'_{1}) +
  \mbox{\rm tr}\,(A_{2}^{t} A'_{2}) +
  \mbox{\rm tr}\,(B^{t}B')
\end{displaymath}
Then: $N_{(A_{1}, A_{2}, B)}{\cal O}(A_{1}, A_{2}, B)\cap {\cal X}$ is
the vector subspace consisting of triples $(X_{1}, X_{2}, Y)\in
{\cal X}$ such that
\begin{displaymath}
 A_{1}X_{1}^{t}-X_{1}^{t}A_{1}+A_{2}X_{2}^{t}-X_{2}^{t}A_{2}-B Y^{t}\in
{\cal A}
\end{displaymath}
where ${\cal A}$ is the set
\begin{displaymath}
\left \{M=(m_{i}^{j})\, \left \vert \, m_{i}^{j} = 0,\ 2 \leq i\leq n,
      1\leq j\leq n \right .\right \}
\end{displaymath}
\end{prop}

Normal spaces of two equivalent triples can be obtained one from the
other. Thus, it is always possible to restrict ourselves to the case
where the triple is in its canonical form.

A {\it bifurcation diagram} of a family of bimodal systems,
\begin{displaymath} \Lambda:
\mathbb{R}^{d}\longrightarrow M_{n}(\mathbb{R})\times
M_{n}(\mathbb{R})\times M_{n\times 1}(\mathbb{R})
\end{displaymath}
is the partition of the parameter space $\mathbb{R}^{d}$ induced by the
stratification associated to the canonical form of the triples of
matrices (see Section~\ref{sec:BPWLS}). In particular, this
stratification provides the information about which canonical forms are
near each other in the sense of local perturbations. Since small changes
in the coefficients of the matrices defining the system may give rise to
matrices defining non-equivalent systems, it is necessary, in order to
explain the behavior of the system under small perturbations, to know
the nearby equivalence classes. Recall that most generic equivalence
classes correspond to lowest codimension in the closure hierarchy.

Let us show how local bifurcation diagrams can be obtained by means of
the miniversal deformation above.

\begin{ex}\label{ex:canonical} Consider a bimodal linear system of type
  CF10' whose canonical form is
\begin{displaymath}
A_{1} = \left(\begin{array}{rr}
                 a_{4} & 0\\
                    0  & a_{4}
               \end{array}\right),\;
A_{2}=\left(\begin{array}{rr}
                 a_{4} & 0\\
                    1  &a_{4}
            \end{array}\right),\;
B=\left(\begin{array}{c}
                 b_{1}\\
                    0
         \end{array}\right).
\end{displaymath}
Then, $N_{(A_{1}, A_{2},B)}{\cal O}(A_{1}, A_{2}, B)\cap {\cal X}$ is
the vector subspace consisting of triples $(X_{1}, X_{2}, Y)\in {\cal X}$
\begin{displaymath}
   X_{1}=\left(\begin{array}{rr}
                  x_{1} & x_{3}\\
                  x_{2} & x_{4}
               \end{array}\right),\;
   X_{2}=\left(\begin{array}{rr}
                  x_{5} & x_{3}\\
                  x_{6} & x_{4}
                 \end{array}\right),\;
   Y=\left(\begin{array}{c}
                  y_{1}\\
                  y_{2}\end{array}\right)
\end{displaymath}
such that
\begin{displaymath}
\left.
\begin{array}{rl}
x_6&=0\\
a_4x_5+b_1y_2&=0\end{array}
\right\} 
\end{displaymath}


Moreover, parameter $x_{3}$ must be zero to avoid observable
perturbations and parameters $x_4,y_1$ give orbits in the initial
stratum.

Then the unobservable perturbations in the normal space to the stratum
of $(A_{1}, A_{2}, B)$ are parameterized by
\begin{displaymath}
  \varphi(x_{1}, x_{2}, x_{5}) = \left(\left(
      \begin{array}{cc}
         a_{4} + x_{1} & 0\\
         x_{2} & a_{4}\end{array}\right),
   \left(\begin{array}{cc}
         a_{4} + x_{5} & 0\\
         1 & a_{4}\end{array}\right ),
   \left (\begin{array}{c}
         b_{1}\\
         -\frac{a_{4}}{b_{1}}x_{5}\end{array}\right)\right)
\end{displaymath}
We denote by $E_{i}$ the set of all triples of matrices having canonical
form of type (CFi), $i=1,\ldots,14$.

Clearly, if only $x_{1}$ (respectively $x_{2}$) is non-zero, it lies in E5'
(respectively E8). But for only $x_{5}$, the strata E6 and E7 are possible
in principle, depending on the value of $\Delta_{2}$. In our case
\begin{displaymath}
\Delta_{2} = b_{1}\gamma_{2} + (a_{4}-\gamma_{1})b_{2} = b_{1} +
 (-x_5)(-\frac{a_{4}}{b_{1}}x_{5}) = \frac{1}{b_{1}}(b_{1}^{2} +
  a_{4}x_{5}^{2}).
\end{displaymath}
Hence, it belongs to E7 for $x_{5}^{2}=-b_{1}^{2}/a_{4}$, and to E6
otherwise.

In a similar way, if $x_{1}, x_{5}\neq 0$ only E2, E3 and E4 are
possible. We have $\Delta_{0} = -x_{2}x_{5} + x_{1}$. Hence, $x_{2}=0$ implies
$\Delta_{0}\neq 0$, which corresponds to E2. If $x_{2}\neq 0$, it gives
again E2 except on the hyperbolic paraboloid $x_{1}=x_{2}x_{5}$. When it
happens:
\begin{displaymath}
\Delta_{1}=b_{1}x_{2}+x_{1}\frac{a_{4}}{b_{1}}x_{5} =
\frac{x_{2}}{b_{1}}(b_{1}^{2}+a_{4}x_{5}^{2}).
\end{displaymath}
Hence, it lies in E4 for $x_{5}^{2}=-b_{1}^{2}/a_4$, and in E3 otherwise.

Finally, it is straightforward that one obtains E5 for $x_1=0$,
$x_2,x_5\neq 0$, and E5' for $x_5=0$, $x_1,x_2\neq 0$.

Summarizing (see Fig. 1):

\begin{figure}[!h]
\centering
\includegraphics[scale=0.35]{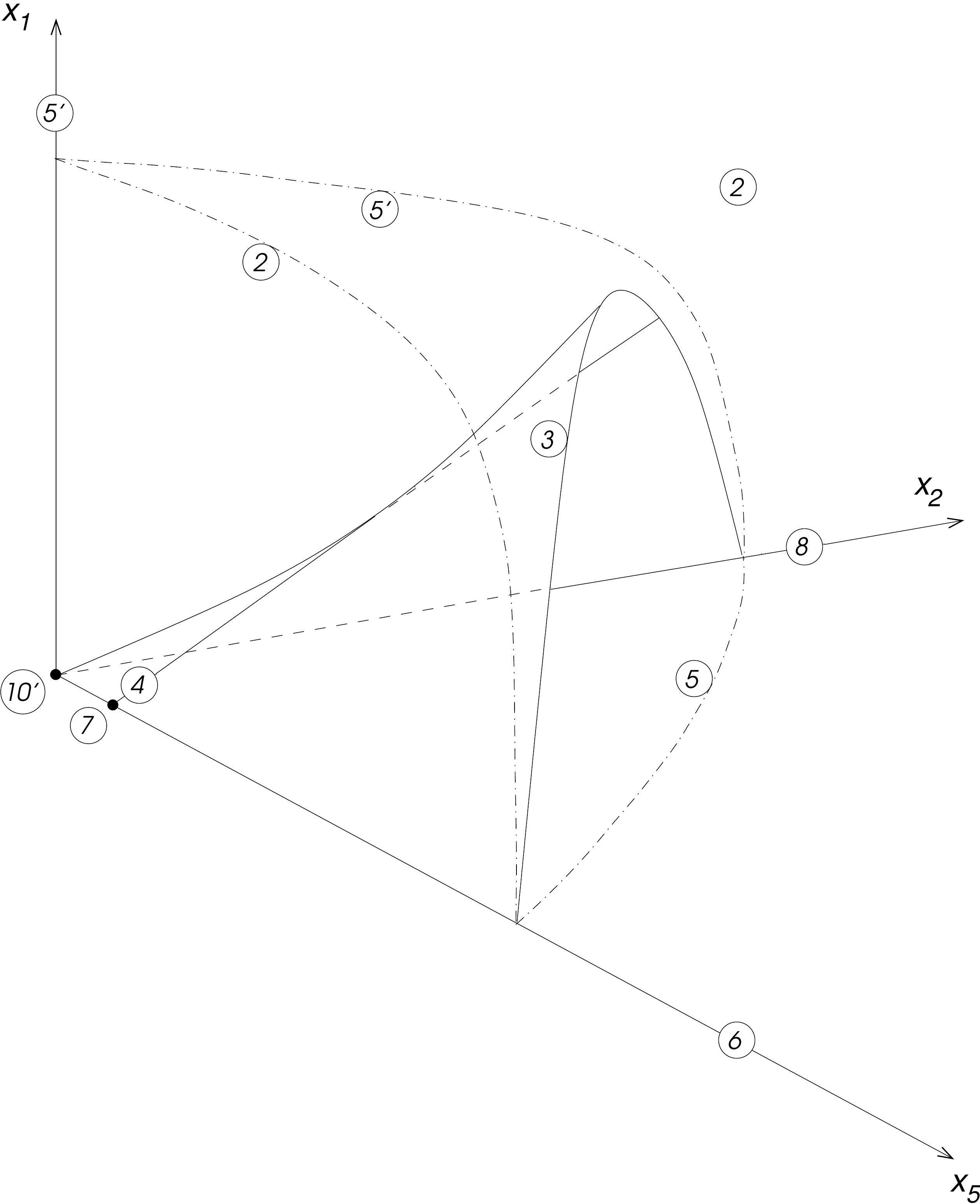}
\caption{Bifurcation diagram.}
\end{figure}


 \begin{enumerate}
\item[$\bullet$] If $x_{2}, x_{5}=0$, $x_{1}\neq 0$, then
  $\varphi(x_{1}, x_{2}, x_{5})\in E5'$.
\item[$\bullet$] If $x_{1}, x_{5} =0$, $x_{2}\neq 0$, then
  $\varphi(x_{1}, x_{2}, x_{5})\in E8$.
\item[$\bullet$] If $x_{1}, x_{2}=0$, $x_{5}^{2} = -b_{1}^{2}/a_{4}$,
  then $\varphi(x_{1}, x_{2}, x_{5})\in E7$.
\item[$\bullet$] If $x_{1}, x_{2} = 0$, $x_{5}\neq 0$, $x_{5}^{2}\neq
  -b_{1}^{2}/a_{4}$, then $\varphi(x_{1}, x_{2}, x_{5})\in E6$.
\item[$\bullet$] If $x_{5}=0$, $x_{1}, x_{2}\neq 0$, then
  $\varphi(x_{1}, x_{2}, x_{5})\in E5'$.
\item[$\bullet$] If $x_{2} = 0$, $x_{1}, x_{5}\neq 0$, then
  $\varphi(x_{1},x_{2},x_{5})\in E2$.
\item[$\bullet$] If $x_{1} = 0$, $x_{2}, x_{5}\neq 0$, then
  $\varphi(x_{1}, x_{2}, x_{5})\in E5$.
\item[$\bullet$] If $x_{1}, x_{2}, x_{5}\neq 0$, $x_{1} = x_{2} x_{5}$,
  $x_{5}^{2} = -b_{1}^{2}/a_{4}$, then $\varphi(x_{1},x_{2},x_{5})\in E4$.
\item[$\bullet$] If $x_{1}, x_{2}, x_{5}\neq 0$, $x_{1} = x_{2} x_{5}$,
  $x_{5}^{2}\neq -b_{1}^{2}/a_{4}$, then $\varphi(x_{1}, x_{2}, x_{5})\in E3$.
\item[$\bullet$] If $x_{1}, x_{2}, x_{5}\neq 0$, $x_{1}\neq x_{2}
  x_{5}$, then $\varphi(x_{1}, x_{2}, x_{5})\in E2$.
 \end{enumerate}


 \end{ex}


\section{Controllability}\label{sec:contr}
As it is well-known, controllability is a qualitative property playing a
central role in many problems. In \cite{CamlibelHS}, one obtains an
implicit characterization of the controllability of bimodal linear
systems. Here, we will characterize explicitly the controllable
unobservable bimodal linear systems for $n=2$ in a quite simple way
(Corollary \ref{cor:contro4}). The canonical forms for bimodal linear
systems in Section \ref{sec:BPWLS} enable simple expressions for these
conditions since they are invariant under admissible basis change
transformations.

\begin{thm} (\cite{CamlibelHS})\label{teor:contro}
Let us consider a bimodal linear system defined by
$(A_{1}, A_{2}, B)$. Let us denote by $e\in M_{n\times 1}(\mathbb{R})$
the matrix such that $eC=A_2-A_1$. Then this system is controllable if,
and only if,
\begin{enumerate}[(1)]
   \item\label{item:contcond} $(A_1,[B\vert e])$ is controllable.
   \item\label{item:vtmui} $\left (v^{t}\ \mu_{i}\right)
          \left(\begin{array}{cc}
                   \lambda I_{n} - A_{i} & B\\
                                      C  & 0\end{array}\right) = 0, \
          \lambda\in \mathbb{R}, v\not= 0$ for $i=1,2\ \Rightarrow \
          \mu_{1}\mu_{2}>0$ holds.
\end{enumerate}
\end{thm}
Next proposition proves that these conditions are invariant under
admissible basis changes, so that they can be checked in the canonical
form of the given system.
\begin{prop}
   Let us consider a controllable bimodal linear system defined by
   $(A_{1}, A_{2}, B)$. Then for all $S\in {\cal S}$, the system
   $(S^{-1}A_{1}S, S^{-1}A_{2}S, S^{-1}B)$ is controllable.
\end{prop}

\noindent {\bf Proof.} First, recall (see Proposition~\ref{prop:pres}
(\ref{item:matrixC})) that the basis change $S$ preserves $C$ (that is,
$C S = C$).

By a standard argument, from condition~(\ref{item:contcond}) it follows
that the system $(S^{-1}AS, [S^{-1}B|S^{-1}e])$ is controllable. Then,
it is sufficient to check
\begin{displaymath}
(S^{-1}e)C=(S^{-1}e)(CS) =
S^{-1}(A_{2}-A_{1})S=S^{-1}A_{2}S-S^{-1}A_{1}S.
\end{displaymath}

Concerning~(\ref{item:vtmui}), let us see that the values $\mu_i$ are
preserved when $S$ acts. Taking into account that
\begin{displaymath}
   \left(\begin{array}{cc}
            \lambda I_{n}-S^{-1}A_{i} S & S^{-1}B\\
                                     CS & 0\end{array}\right) =
   \left(\begin{array}{cc}
            S^{-1} & \\
                   & I_n\end{array}\right)
   \left(\begin{array}{cc}
            \lambda I_{n}-A_{i} & B\\
                             C  & 0\end{array}\right)
   \left(\begin{array}{cc}
            S & \\
            & I_{n}\end{array}\right)
\end{displaymath}
it is obvious that
\begin{displaymath}
   (v^t \ \mu_i) \left(\begin{array}{cc}
                          \lambda I_{n}-S^{-1}A_{i}S & S^{-1}B\\
                                                  CS &0\end{array}
                 \right) = 0
\end{displaymath}
is equivalent to
\begin{displaymath}
   (v^{t}S^{-1} \ \mu_{i})\left(\begin{array}{cc}
                                   \lambda I_{n}-A_{i} & B\\
                                                    C  & 0\end{array}
                          \right) = 0
\end{displaymath}
\hfill $\Box$


Let us assume that $n=2$ and consider the unobservable system defined by
$(A_{1}, A_{2}, B)$, where
\begin{displaymath}
A_{1} = \left(\begin{array}{cc}
                 a_{1} & 0\\
                 a_{2} & a_{4}\end{array}\right),
A_{2} =\left(\begin{array}{cc}
                 \gamma_{1} & 0\\
                 \gamma_{2} & a_{4}\end{array}\right),
B=\left(\begin{array}{c}
           b_{1}\\
           b_{2}\end{array}\right)
\end{displaymath}
The following theorem gives an equivalent condition
to~(\ref{item:vtmui}) in Theorem~\ref{teor:contro}:

\begin{prop}\label{prop:contro2}
Let us consider an unobservable bimodal linear system defined by
$(A_{1}, A_{2}, B)$. With the above notation, if $n=2$,
condition~(\ref{item:vtmui}) in Theorem~\ref{teor:contro} is equivalent
to
\begin{displaymath}
   b_{1}\neq 0,\quad \Delta_{1}\Delta_{2} > 0.
\end{displaymath}
\end{prop}

\noindent {\bf Proof.}
Condition~(\ref{item:vtmui}) in Theorem~\ref{teor:contro} may be
re-written as follows:
%
\begin{displaymath}
\left .
\begin{array}{rl}
b_{1} v_{1} + b_{2} v_{2} & = 0\\
v_{2} (\lambda - a_{4}) & = 0\\
v_{1} (\lambda - a_{1}) - v_{2} a_{2} + \mu_{1} & = 0
\end{array}\right\}
\end{displaymath}

\begin{displaymath}
\left.
\begin{array}{rl}
b_{1} v_{1} + b_{2} v_{2} & = 0\\
v_{2}(\lambda - a_{4}) & = 0\\
v_{1}(\lambda-\gamma_{1}) - v_{2} \gamma_{2} + \mu_{2} & = 0
\end{array}\right\}
\end{displaymath}
$\lambda\in\mathbb{R}$, $(v_{1}, v_{2})\neq (0, 0)$ implies
$\mu_{1}\mu_{2} > 0$.

\begin{enumerate}[(a)]
\item If $b_{1}\neq 0$, then
\begin{displaymath}
   v_{1} = -b_{2}, \quad v_{2} = b_{1}, \quad \lambda = a_{4}.
\end{displaymath}
Hence,
\begin{displaymath}
   \mu_{1} = \Delta_{1},\quad \mu_{2} = \Delta_{2}
\end{displaymath}
%
\item If $b_{1} = 0$,
%
%
the first and the second equations are verified by
\begin{displaymath}
       v_{1}\neq 0,\quad v_{2} = 0,\quad
            \text{ any }\lambda\in\mathbb{R}
\end{displaymath}
so that any $\mu_{i}$ is possible, and~(\ref{item:vtmui}) does not
hold.
\end{enumerate}
\hfill $\Box$

\begin{cor}\label{cor:contro3}
In the conditions of Theorem~\ref{teor:contro}, if $n=2$,
condition~(\ref{item:vtmui}) implies condition~(\ref{item:contcond}).
\end{cor}

\noindent {\bf Proof.}
Condition~(\ref{item:contcond}) in Theorem~\ref{teor:contro} may be
re-written as follows:
\begin{displaymath}
\mbox{\rm rk}\,\left(\begin{array}{cccc}
                        b_{1} & a_{1}b_{1} & \gamma_{1}-a_{1} &
                        a_{1}(\gamma_{1}-a_{1})\\
                        b_{2} & a_{2}b_{1} + a_{4} b_{2} &
                              \gamma_{2} - a_{2} &
                              a_{2}(\gamma_{1} - a_{1}) +
                              a_4(\gamma_2-a_2)
                     \end{array}\right)=2,
\end{displaymath}
which clearly follows from $b_{1}\Delta_{1}\neq 0$.\hfill $\Box$
%
%
%

Our first main result follows from Theorem~\ref{teor:contro},
Proposition~\ref{prop:contro2} and Corollary~\ref{cor:contro3}.

\begin{cor}\label{cor:contro4}
Let us consider an unobservable bimodal linear system defined by
$(A_1,A_2,B)$. If $n=2$, this system is controllable if, and only if,
$$b_1\neq 0,\quad \Delta_1\Delta_2>0.$$
\end{cor}
%
%

{\bf Remark}
From this formulation and Lemma \ref{lem:lemn1}~(\ref{item:det}), it is
obvious that if $(A_{1}, A_{2}, B)$ is controllable, then both
subsystems $(A_{1}, B)$ and $(A_{2}, B)$ are controllable as well.
Next Table summarizes the results when the above condition is applied to
the canonical form of each stratum:

\vglue.5truecm

\centerline{\kern-1.2truecm
\begin{tabular}{|c|c|}
  \hline
Canonical form&Controllability\\\hline
CF2 & $b_{1}\ne 0$ and $\Delta_{1}\Delta_{2} > 0$\\\hline
CF3 & $b_{1}\ne 0$ and $(a_{4}-a_{1})(a_{4}-\gamma_{1}) > 0$\\\hline
CF4 & Always uncontrollable \\ \hline
CF5 & $b_{1}\Delta_{2} > 0$\\ \hline
CF5'& $b_{1}\Delta_{1} > 0$\\ \hline
CF6 & Always uncontrollable \\ \hline
CF6'& Always uncontrollable \\ \hline
CF7 & Always uncontrollable \\ \hline
CF7'& Always uncontrollable \\ \hline
CF8 & $a_{2}\gamma_{2} > 0$ \\ \hline
CF9 & Always uncontrollable \\ \hline
CF10 & Always uncontrollable \\ \hline
CF10' & Always uncontrollable \\\hline
CF11 & Always uncontrollable \\ \hline
CF11' &Always uncontrollable \\ \hline
CF12 & Always uncontrollable \\ \hline
CF13 & Always uncontrollable \\ \hline
CF14 & Always uncontrollable\\\hline
\end{tabular}
}

Let us show an example to illustrate the study of the controllability of
a bimodal linear system.

\begin{ex}\label{ex:canonical2} Consider a bimodal linear system of type
  CF10' whose canonical form is
\begin{displaymath}
   A_{1}=\left(\begin{array}{rr}
                  a_{4} & 0\\
                      0 & a_{4}\end{array}\right),\;
   A_{2}=\left(\begin{array}{rr}
                  a_{4} & 0\\
                      1 & a_{4}\end{array}\right),\;
   B=\left(\begin{array}{c}
                  b_{1}\\0\end{array}\right).
\end{displaymath}

We consider the unobservable perturbation obtained in
Example~\ref{ex:canonical}:
\begin{displaymath}
   \varphi(x_{1}, x_{2}, x_{5}) =
      \left(\left(\begin{array}{cc}
                     a_{4} + x_{1} & 0\\
                             x_{2} & a_{4}\end{array}\right),
      \left(\begin{array}{cc}
                     a_{4} + x_{5} & 0\\
                                 1 & a_{4}\end{array}\right),
      \left(\begin{array}{c}
                     b_{1}\\
                     -\frac{a_{4}}{b_{1}} x_{5}\end{array}\right )
      \right )
\end{displaymath}
The controllable bimodal linear systems are those satisfying the
condition in Corollary \ref{cor:contro4}
which, taking into account that $b_{1}\neq 0$, is equivalent to
\begin{displaymath}
   \left(x_{2} + \frac{a_{4}}{b_{1}^{2}}x_{1} x_{5}\right)
   \left(1+\frac{a_{4}}{b_{1}^{2}}x_{5}^{2}\right) >0,
\end{displaymath}
which can be decomposed into:
\begin{itemize}
\item If $a_{4} > 0$: \quad $x_{2} > -
                      \frac{a_{4}}{b_{1}^{2}} x_{1} x_{5}$
\item If $a_{4} < 0$: \quad $x_{2} > -
                      \frac{a_{4}}{b_{1}^{2}} x_{1} x_{5}$ and
                      $ x_{5}^{2}> -\frac{b_{1}^{2}}{a_4}$,\quad or
                      \quad $x_{2} < -\frac{a_{4}}{b_{1}^{2}}
                               x_{1} x_{5}$ and $x_{5}^{2}
                               < -\frac{b_{1}^{2}}{a_4}$
\end{itemize}

If $a_{4} < 0$, the above condition is illustrated in
Figure~\ref{fig:controllability}. Figure~\ref{fig:canonical12} summarizes
examples~\ref{ex:canonical} and~\ref{ex:canonical2}.
%
\begin{figure}[!t]
\centering
\includegraphics[scale=0.6]{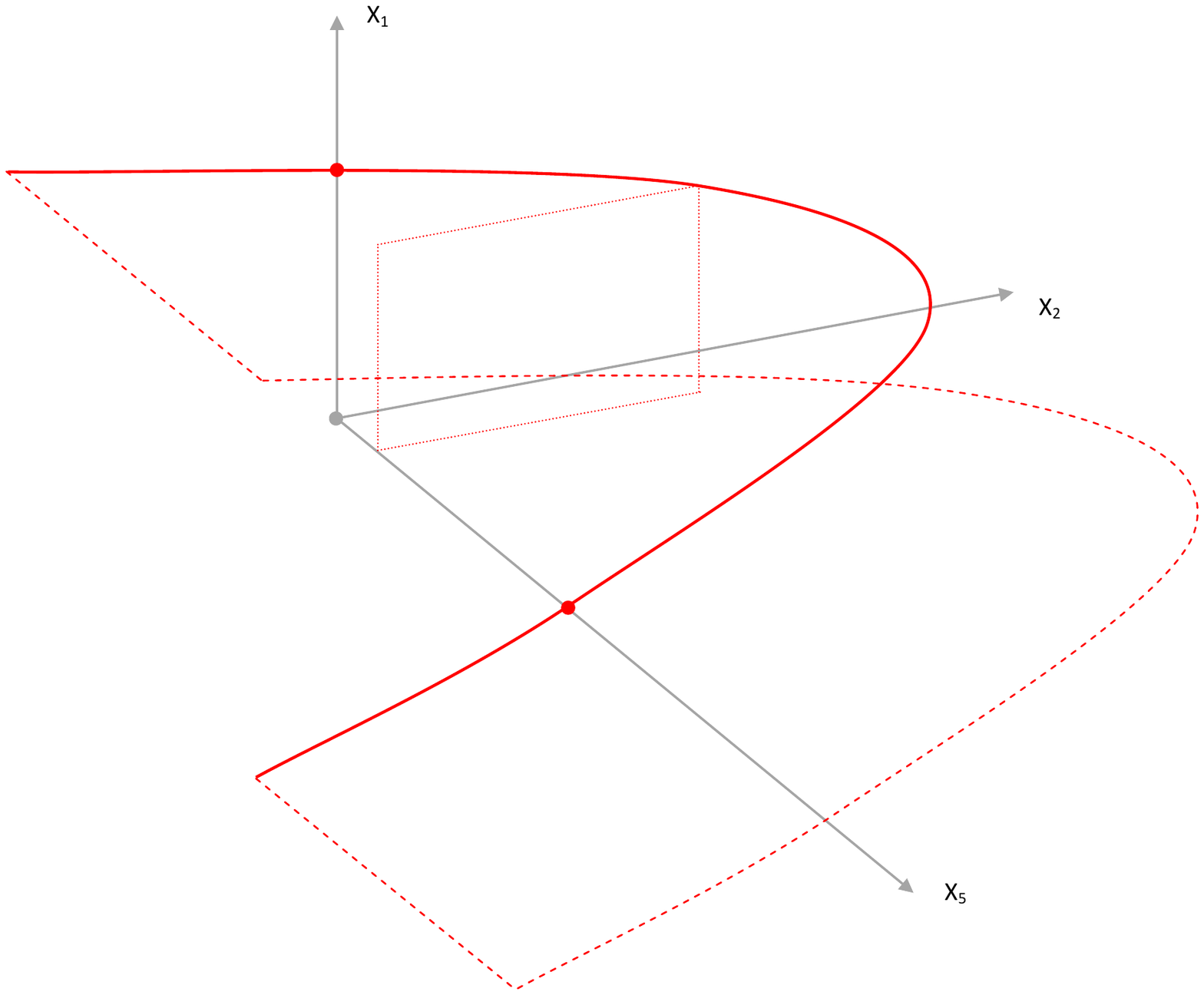}
\vspace{-15\baselineskip}
\caption{Controllability.\label{fig:controllability}}
\end{figure}



%

\begin{figure}[!t]
\centering
\includegraphics[scale=0.6]{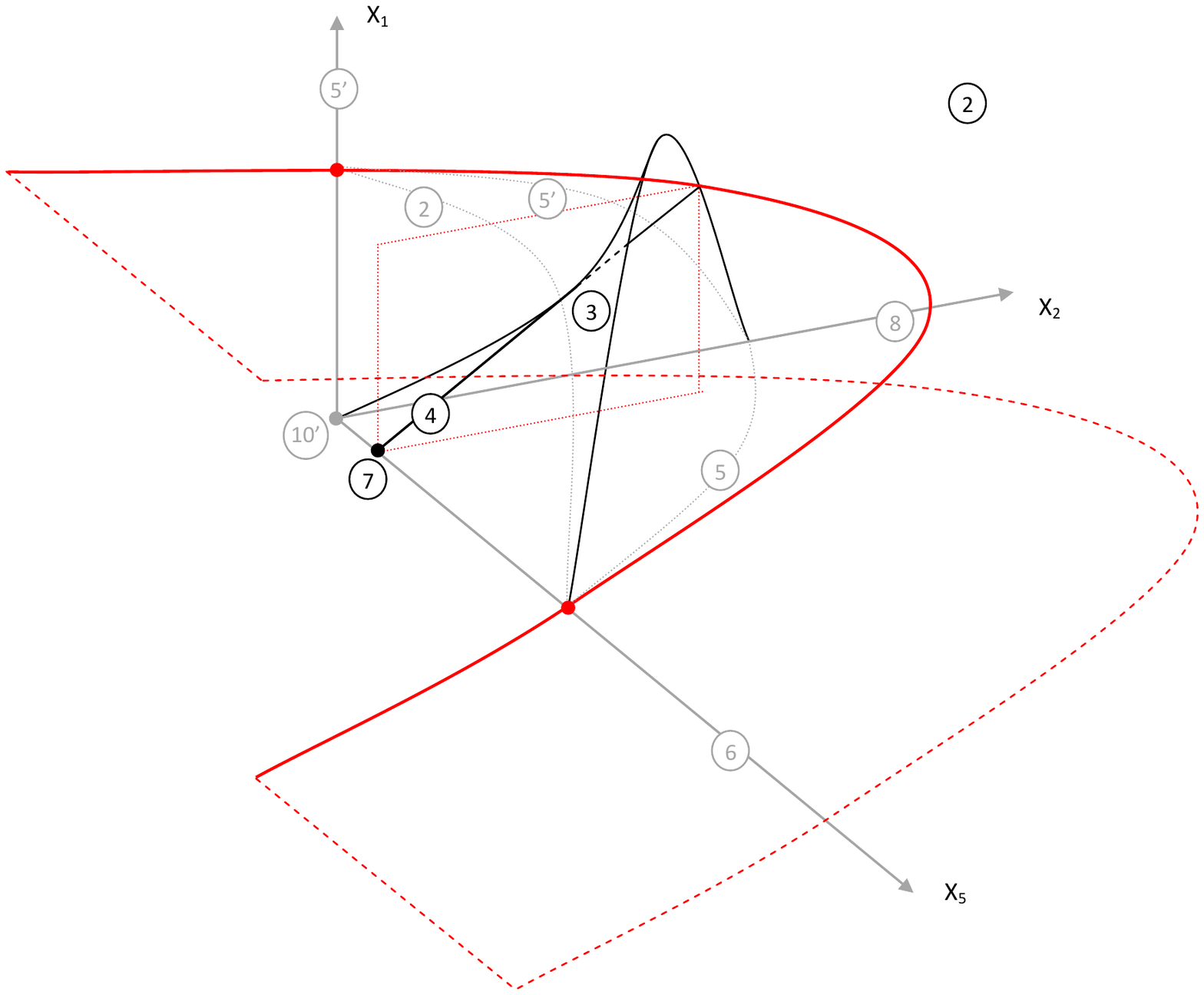}
\vspace{-15\baselineskip}
\caption{Bifurcation diagram and Controllability.~\label{fig:canonical12}}
\end{figure}
\end{ex}
Finally, we will use the characterization of controllability in
Corollary \ref{cor:contro4} and the canonical forms in
Section~\ref{sec:BPWLS} to prove that then the system is stabilizable,
by means of a common feedback for both subsystems.


\begin{thm}\label{teor:stabi}
Let us consider an unobservable planar bimodal linear system defined by
$(A_{1}, A_{2}, B)$. If it is controllable, there is a feedback $F\in
M_{1\times 2}(\mathbb{R})$ such that both subsystems $A_{1} + BF$ and
$A_{2} + BF$ are stable.
\end{thm}

\noindent
{\bf Proof.}
We will detail the proof for the canonical form CF2. It works
analogously for CF3, CF5, CF5' and CF8.

By Corollary~\ref{cor:contro4}, the system
\begin{displaymath}
A_{1} = \left(\begin{array}{rr}
                 a_{1} & 0\\
                     0 & a_{4}\end{array}\right),\;
A_{2} = \left(\begin{array}{rr}
                 \gamma_{1} & 0\\
                          1 & a_{4}\end{array}\right),\;
B = \left(\begin{array}{c}
             b_{1}\\
             b_{2}\end{array}\right)
\end{displaymath}
is controllable if and only if
\begin{enumerate}[(i)]
\item\label{item:b1ne0} $b_{1}\neq 0$
\item\label{item:Delta1} $\Delta_1\Delta_2>0$
\end{enumerate}
where $\Delta_{1} = (a_{4}-a_{1})b_{2},\quad \Delta_{2} =
       b_{1} + (a_{4} - \gamma_{1})b_{2}$.

We search $F=\left(\begin{array}{cc} f_{1} & f_{2}\end{array}\right)$
such that the matrices $A_{1}+BF$ and $A_{2}+BF$ have negative trace and
positive determinant, that is to say:
\begin{gather*}
b_{1} f_{1} + b_{2}f_{2} < -a_{1} - a_{4}\\
b_{1} f_{1} + b_{2}f_{2} < -\gamma_{1} -a_{4}\\
a_{4}b_{1}f_{1} + a_{1}b_{2}f_{2} > -a_{1}a_{4}\\
a_{4}b_{1}f_{1} + \gamma_{1}b_{2}f_{2}-b_{1}f_{2}> - \gamma_{1}a_{4}
\end{gather*}
We can change the variables $(f_{1}, f_{2})$ by $(x, y)$ defined by
\begin{eqnarray*}
x &=& b_1f_1+b_2f_2\\
y &=& -(a_4b_1f_1+a_1b_2f_2)
\end{eqnarray*}
because (recall~(\ref{item:b1ne0}) and~(\ref{item:Delta1}))
\begin{displaymath}
   \det\left(\begin{array}{cc}
     b_{1} & b_{2}\\
    -a_{4} b_{1} & -a_{1}b_{2}\end{array}\right)= b_{1}\Delta_{1}\neq 0
\end{displaymath}
Then:
\begin{displaymath}
  f_{1} = -\frac{a_{1}x + y}{(a_{4} - a_{1})b_{1}},\quad
  f_{2} = \frac{a_{4}x + y}{(a_{4} - a_{1})b_{2}}
\end{displaymath}
With this change of variables, the desired inequalities become:
\begin{gather*}
 x < -a_{1} -a_{4}\\
 x < -\gamma_{1} - a_{4}\\
 y < a_{1}a_{4}\\
 a_{4}b_{1}\frac{a_{1}x + y}{(a_{4} - a_{1})b_{1}} +
      (b_{1} -\gamma_{1}b_{2})\frac{a_{4}x + y}{(a_{4} - a_{1})b_{2}} <
      \gamma_{1}\gamma_{4}
\end{gather*}
In order to see that there exist solutions $(x,y)$, it is sufficient
that the coefficient of the variable $y$ in the last inequality be
positive:
\begin{displaymath}
   \frac{a_{4} b_{1}}{(a_{4} - b_{1})b_{1}} +
   \frac{b_{1} - \gamma_{1}b_{2}}{(a_{4}-a_{1})b_{2}} =
   \frac{1}{b_{1}\Delta_{1}}(a_{4}b_{1}b_{2}+b_{1}^{2} -
   \gamma_{1}b_{1}b_{2}) =
   \frac{\Delta_{2}}{\Delta_{1}} > 0
\end{displaymath}
again by~(\ref{item:Delta1}).\hfill $\Box$

\subsection*{Acknowledgements}
Preprint of an article submitted for consideration in IJBC \copyright
2011 copyright World Scientific Publishing Company
\texttt{http://www.worldscinet.com/ijbc/}


\begin{thebibliography}{1}
\bibitem{Arnold71} V.~I.~Arnold, \newblock{On matrices depending on
    parameters}. \newblock Uspekhi Mat.~Nauk., 26
(1971).
\bibitem{CamlibelHS} K.~Camlibel, M.~Heemels, H.~Schumacher, On the
  controllability of bimodal piecewise linear systems, LNCS 2993 (2004),
  p. 250--264.
\bibitem{CarmonaFPT02} V.~Carmona, E.~Freire, E.~Ponce,
F.~Torres, \newblock{On simplifying and classifying piecewise linear
  systems}. \newblock IEEE Transactions on Circuits and Systems, 49
(2002), p.~609--620.
\bibitem{CarmonaFPT06} V.~Carmona, E.~ Freire, E.~Ponce, F.~Torres,
\newblock{The con\-ti\-nuous matching of two stable linear systems can
           be unstable}.
\newblock Discrete and continuous dynamical systems, 16, 3 (2006),
          p.~689--703.
\bibitem{DiBernardoPP08} M.~di Bernardo, D.~J.~Pagano, E.~Ponce,
\newblock{Nonhyperbolic boundary equili\-brium bifurcations in planar
  {F}ilippov systems: a case study approach}. \newblock In\-ter\-nat.
J.~Bifur.~Chaos Appl.~Sci.~Engin., 18, 5 (2008), p.~1377--1392.
\bibitem{DiBernardoBCK08} M.~di Bernardo, C.~J.~Budd, A.~Champneys,
                     P.~Kowalczyk,
\newblock{\em Piecewise-Smooth Dynamical Systems}. \newblock
Springer-Verlag, London (2008).
 \bibitem{MP} J.~Ferrer, M.~D.~Magret, M.~Pe\~ na,
\newblock{Bi\-mo\-dal pie\-ce\-wi\-se li\-near systems. Reduced
  Forms}. \newblock International Journal of Bifurcation and Chaos, 20,
9 (2010),  p.~2795--2808.
\bibitem{MPP} J.~Ferrer, M.~D.~Magret, J.\,R. Pacha,
  M.~Pe\~na,\newblock{Planar Bimodal Piecewise Linear
    Systems. Bifurcation  Diagrams}. Bol. Soc. Esp. Mat. Apl., 51
  (2010), p.~55--63.
\bibitem{FreirePR05} E.~Freire, E.~Ponce, J.~Ros, \newblock{The
    focus-center-limit cycle bi\-fur\-ca\-tion in sym\-me\-tric 3D
    piecewise linear systems}.
\newblock SIAM J. Appl.~Math., 65, 3 (2005), p.~1933--1951.
\bibitem{Hum}
J.~E.~Humphreys, \newblock{\em Linear Algebraic Groups}. \newblock
Graduate Texts in Mathematics, 21, Springer-Verlag, Berlin (1981).
\bibitem{tan} {A. Tannenbaum},
{\em Invariance and system theory: algebraic and geometric
aspects}, {LNM, n. 845, Springer Verlag} (1981).
 \end{thebibliography}
\end{document}